## There are natural scores: Full comment on Shafer,
## 'Testing by betting: A strategy for statistical and scientific communication'

Sander Greenland, Professor Emeritus

Department of Epidemiology and Department of Statistics

University of California, Los Angeles, California, U.S.A.

10 February 2021

I promote a perspectivist (toolkit) view of statistics, in which statistical philosophies are treated as ways of looking at data, and the methods they lead to are then treated as tools whose utility is determined on an application-specific basis. Shafer offers a betting perspective on statistical testing which may be useful for foundational debates, given that disputes over such testing continue to be intense (e.g., see the special issues of TAS introduced by Wasserstein and Lazar 2016 and Wasserstein et al. 2019). To be helpful for researchers, however, I believe the perspective will need more elaboration using real examples in which (a) the betting score has a justification and interpretation in terms of study goals that distinguishes it from the uncountable mathematical possibilities, and (b) the assumptions in the sampling model are uncertain.

On justification, Shafer says "No one has made a convincing case for any particular choice" of a score derived from a P-value and then states that "the choice is fundamentally arbitrary." I think that is quite wrong in that some (but not most) scores can be motivated by study goals (e.g., information measurement; decision making). The one choice Shafer displays $(p^{-\frac{1}{2}} - 1)$ does appear to have no justification, but I find it most unnatural; its absence of history suggests others do as well. Yet there are justifiable choices. The one I have seen repeatedly in information statistics and data mining is the surprisal, logworth or S-value $s_b = \log_b(1/p) = -\log_b(p)$, where the log base b determines the scale via $s_b = s_e/\ln(b)$ (Fisher 1948; Good 1956, 1957; Bayarri and Berger 1999; Boos and Stefanski 2011; Greenland 2019; Fraundorf 2020).

To explain this choice, consider one traditional version statistical testing of a simple hypothesis H, which assumes

1) The available information about the mechanism that generated the data can be represented adequately by a set of background assumptions about the data generator (e.g., uncorrelated residuals) which define a model space A.

2) The scientific hypothesis under study can be represented adequately by a set of constraints H on A, called a test hypothesis, which defines a subset HA of A.

We then find a measure D of the discrepancy between the data sets expected under A and under HA that is 'tractable': the distribution $F_{HA}(d)$ of D under HA is known or identifiable. Typically D is a geometric or information distance of the projection of the data onto A from its projection onto HA. Then $P = 1-F_{HA}(D)$ is the (random) upper-tail P-value and is uniform when HA holds, and the observed 'exact' P-value is the percentile $p = 1-F_{HA}(d)$. It follows that P is uniform when HA holds; the decision rule "reject H if $p \leq \alpha$" has a rejection rate (size) of $\alpha$; and $S_e = -\ln(P)$ has a unit-exponential distribution with expectation 1. $S_e$ is a betting score in Shafer's system, the "natural" score derived from p (in the narrow sense of the log base).

Using instead $b = 2$ we scale up the score to $s_2 = -\log_2(p) = s_e/\ln(2)$, which is the Shannon information (in bits) or surprisal index for the event D≥d under HA. Inverting back to $p = 2^{-s_2}$ provides an intuitive calibration in terms of fair coin tossing: Observations of p = 0.10, 0.05, 0.01, 0.005, 0.001, 0.000001 become $s_2 = -\log_2(p) = s_e/\ln(2) = 3.32, 4.32, 6.64, 7.64, 9.97, 19.9$; after rounding each $s_2$, we see that each p corresponds roughly to the probability of seeing all heads in 3, 4, 7, 8, 10, or 20 such tosses. More generally, $s_2$ can be treated as measuring the number of bits of Shannon information provided by D≥d against H given A, or against HA unconditionally (Greenland 2019; Rafi and Greenland 2020; Fraundorf 2020). For conservative tests, $s_2$ becomes the minimum information against HA provided by D≥d.

Bayarri and Berger (1999) and Sellke et al. (2001) provide other calibrations based on certain restrictions, e.g., a lower bound of $s_e/\exp(s_e-1)$ for a Bayes factor, subject to certain restrictions on the prior distribution. Direct calibrations of P-values and S-values against sampling experiments have an advantage over Bayesian formulations in not requiring such restrictions. Still, the utility of any test is largely determined by its sensitivity to departures from both H *and* the background assumptions in A used to compute p. Sensitivity to departures from H is usually measured by power functions of $\alpha$-level decision rules under A; while robustness is usually measured by *in*sensitivity of test size and power to departures from assumptions in A. More generally, test sensitivity can be measured by shifts in the entire P or S distribution under violations of H given A (e.g., Sellke et al. 2001; Boos and Stefanski 2011), violations of A given H (when H remains defined when A is violated), or violations of both.

On interpretation and terminology, the observed P-value p has long been used as an index of consistency or *compatibility* of H with data, given A (e.g., Bayarri et al. 2000; Robins et al. 2000), while reversing transforms the such as $1-p$, $1/p$, and $s_b$ have been used as measures of evidence against H. These uses are reasonable if we are certain of A, as we might be if we conducted an experiment in which every condition in A was successfully enforced by study conduct and laws of nature. But if we are unwilling to condition our inferences on A (e.g., because we are concerned about possibly important protocol violations), we cannot say more than that they measure evidence against HA. In that case, the Shannon information $s_2$ is only measuring the bits of information supplied by the test against HA; a large value alone does not allow us to say which of H or A is violated, as it may arise from violations of H, A, or both (Greenland 2019; Greenland and Rafi 2020). Other tests may be more specific, but those will require assumptions of their own. This critical, *unconditional* view of statistics extends to Bayesian settings by including the prior distributions in A (as random-parameter distributions in the data-generating model); in that case the information will now be against the model for the joint distribution of the data and the model parameters (Box 1980).

I believe the unconditional view is far better suited to observational studies than are traditional statistical descriptions (Greenland and Rafi 2021). The latter focus on H and take too much of A for granted, even though A is often filled with unsupported assumptions to create identification of target parameters (Greenland 2005). For causal targets, one such assumption is that treatment has been assigned randomly given some identifiable propensity score or baseline-risk score (often encoded as "conditional ignorability of treatment assignment" or "no unmeasured confounding"). In contrast, the unconditional description more openly allows for uncertainties about treatment assignment. Parallel comments apply to nonrandom survey selection and nonrandomly missing data, where target identification involves uncertain assumptions of random selection or missingness.

In closing, I agree with Shafer about terminology, insofar as statistics has taken ordinary words like "significance" and "confidence" and turned them into overconfident jargon for methods that neglect crucial information about study problems (Amrhein et al. 2019; Rafi and Greenland 2020; Greenland and Rafi 2021). In the mid-20th century, scientific communities adopted this jargon enthusiastically; yet despite the corruption of the literature it produced, some opinion leaders still defend it vociferously – unsurprising, as to formally account for study problems undermines the significance and confidence that should be assigned to many studies, including their own. Nonetheless, although I sense Shafer may be offering a valuable new viewpoint, I would need to see how Shafer's approach translates into the above logical framework before I could consider using or teaching it, especially in terms of what it means when background assumptions are in doubt (Greenland and Rafi 2021). I suspect others will need similar translations into their own preferred frameworks.